\documentstyle[12pt,amssymb]{article}
\newcommand{\rest}{\restriction}
\newcommand{\dom}{{\mathrm {dom}}}
\newcommand{\supp}{{\mathrm {supp}}}
\newcommand{\QED}{\vrule width 6pt height 6pt depth 0pt \vspace{0.1in}}

\newtheorem{theorem}{Theorem}[section]
\newtheorem{lemma}[theorem]{Lemma}
\newtheorem{corollary}[theorem]{Corollary}
\newtheorem{claim}[theorem]{Claim}
\begin{document}

\noindent Tomek Bartoszynski, Department of Mathematics, Boise State
University, 
Boise, ID 83725
\footnote{Research partially supported by
SBOE grant 92--096}
\bigskip
\begin{center}
{\bf{\LARGE{On the structure of measurable filters on a  countable
set}}}
\end{center}
\bigskip
\begin{abstract} A combinatorial characterization of measurable filters on
a countable set is found. We apply it to the problem of measurability of
the intersection of nonmeasurable filters.
\end{abstract}

The goal of this paper is to characterize measurable filters on the set of
natural numbers. In  section 1 we introduce basic notions, in section 2 we find
a combinatorial characterization of measurable filters, in section 3 we study
intersections of filters and finally section 4 is devoted to filters which are
both null and meager.

Through this paper we use standard notation.
 $\omega$ denotes the set of natural
numbers. For $k,n \in \omega$  let
$[n,k] = \{i<\omega:n \leq i \leq k\}$. For $ n \in \omega , \ 2^{n}\
(2^{\omega})$ denotes the
set of 0-1 sequences of length $n (\omega)$, also let
$2^{<\omega} = \bigcup_{n \in \omega}2^{n}$. For any
sequences $s,t \in 2^{<\omega}$ let $s ^{\frown}t$ denote their
concatenation. For $s \in 2^{<\omega}$ let $[s]=\{x \in 2^{\omega}:
s \subset x\}$. The family $\{[s]:s \in 2^{<\omega}\}$ is a base
of the space $2^{\omega}$. We will often identify
a set $[s]$ with a sequence s and we
will also identify subsets of $\omega$ with their
characteristic functions.
Filters considered in this paper are assumed to be nonprincipal. We identify
filters on $\omega$ with sets of characteristic
functions of its elements. In this way
the question about measurability makes sense.
Finally let quantifiers ``$\exists^{\infty}$'' and ``$\forall^{\infty}$''
denote ``for infinitely many'' and
``for all
except finitely many'' respectively.

\section{Introduction}
In this section we establish some
definitions and recall several facts which
we will use later. Let us start with measures we will be working with.

{\sc Definition}
Let $\hat{p} = \{p_{n}:n \in \omega \}$ be a sequence of
reals such that $p_{n} \in (0,\frac{1}{2}]$ for all $n \in \omega$.
Define $\mu_{\hat{p}}$ to be the product measure on
$2^{\omega}$ such that $\mu_{\hat{p}}(\{x \in 2^{\omega} : x(n)=1\}) =
p_{n}$ and
$\mu_{\hat{p}}(\{x \in 2^{\omega} : x(n)=0\}) = 1 - p_{n}$ for $n \in \omega$.
Notice that if $p_{n} = \frac{1}{2}$ then $\mu_{\hat{p}}$
is the usual measure on $2^{\omega}$.
From now on let us fix one of the measures
$\mu_{\hat{p}}$. We have the following:

\begin{theorem}[Sierpinski]\label{0.1}
Suppose  that ${\cal F}$ is a filter on $\omega$.
Then ${\cal F}$ is either of $\mu_{\hat{p}}$-measure zero or
${\cal F}$ is
$\mu_{\hat{p}}$-nonmeasurable. Moreover, ${\cal F}$
is either meager or does not have the Baire
property.  
\end{theorem}
{\sc Proof}:
In the case of the Baire property or when $\mu_{\hat{p}}$ is the
Lebesgue measure we use the fact that the automorphism of $2^\omega$ which
sends every set to its complement preserves Lebesgue measure.

Suppose that a filter ${\cal F}$ is $\mu_{\hat{p}}$-measurable.
Since ${\cal F}$ is non-principal it has measure 0 or 1.
We have to show that $\mu_{\hat{p}}({\cal F})=0$.

Consider $\varphi : 2^\omega  \times 2^\omega \longrightarrow
2^\omega$ defined as $\varphi(X,Y)(n)=\max(X(n), Y(n))$ for $X,Y \in
2^\omega$.

Let $q_n$ be chosen in such a way that $(1-p_n)(1-q_n) = 1/2$ for all
$n$.

\begin{claim}\label{2.5}
$\varphi^{-1}(\mu) = \mu_{\hat{p}} \times \mu_{\hat{q}}$.
\end{claim}
{\sc Proof}:
Verify that for all $n \in \omega$, 
$$\frac{1}{2}=\mu(\{ x \in 2^{\omega} :
x(n)=0 \}) = \mu_{\hat{p}} \times \mu_{\hat{q}} (\varphi^{-1}(\{x \in
2^{\omega} : x(n) = 0\}))$$ and use the fact that sets of this form are
independent.~$\QED$

Since $\varphi$ is a continuous mapping it follows that if $A
\subseteq 2^\omega \times 2^\omega$ and  $\mu_{\hat{p}}\times
\mu_{\hat{q}}(A)=1$ then $\mu(\varphi(A))=1$.

Consider the set ${\cal F} \times 2^\omega$. If $\mu_{\hat{p}}({\cal
F})=1$ then $\mu_{\hat{p}}\times
\mu_{\hat{q}}({\cal F} \times 2^\omega)=1$. In particular
$\mu(\varphi({\cal F} \times 2^\omega))=\mu({\cal F})=1$. Contradiction.~$\QED$

Sierpinski also proved that if  ${\cal F}$
is an ultrafilter then  ${\cal F}$ is Lebesgue
nonmeasurable. The next theorem shows that with
measures $\mu_{\hat{p}}$ this is not the
case.

\begin{theorem}\label{0.2}
Let ${\cal F}$ be a filter on $\omega$.
\begin{enumerate}
\item If there exists $X \in {\cal F}$ such
that $\sum_{n \in X} p_{n}^{k} < \infty$ for some $k \in \omega$
then $\mu_{\hat{p}}({\cal F}) = 0$.
\item Let $\{k_{n} : n \in \omega \}$ be a sequence of natural numbers
such that $$\sum_{n=1}^{\infty} p_{n}^{k_{n}}~=~\infty
\hbox{ but } \sum_{n=1}^{\infty} p_{n}^{k_{n}+1}<\infty.$$
If for every $X \in {\cal F}$,  $\sum_{n \in X} p_{n}^{k_{n}} = \infty$
then $\mu_{\hat{p}}({\cal F})=0$.
\end{enumerate}
\end{theorem}
{\sc Proof}:
The first part of this theorem is due to M. Talagrand. For completeness we
sketch the proofs of both parts.

1)  If there exists $X \in {\cal F}$ such that
$\sum_{n \in X} p_{n} < \infty$ then
$\mu_{\hat{p}}({\cal F})=0$ since $\mu_{\hat{p}}(\{Y \subseteq \omega:
|X \cap Y|= \omega\})=0$. If there exists $X \in {\cal F}$
such that $\sum_{n \in X} p_{n}^{k+1} < \infty$ but
for all $Y \in {\cal F}$,
$\sum_{n \in Y} p_{n}^{k} = \infty$  then $\mu_{\hat{p}}({\cal F})=0$ since
$\mu_{\hat{p}}(\{Y \subseteq \omega :|Y \cap X|=\omega \hbox{ and }
\sum_{n \in X \cap Y} p_{n}^{k}= \infty \})=0$.

2)   Consider the random variables
$$\xi_{n}(X) = \left\{ \begin{array}{ll}
p_{n}^{k_{n}} & \hbox{if } n \in X \\
0 & \hbox{otherwise}
\end{array}
\right. .$$
Let $\xi = \sum_{n=1}^{\infty} \xi_{n}$. The mean value
$E(\xi) = \sum_{n=1}^{\infty} p_{n}^{k_{n}+1} < \infty$.
Therefore $\mu_{\hat{p}}(\{X \subseteq \omega : \sum_{n \in X}
p_{n}^{k_{n}} = \infty\})=0$.~$\QED$

The next theorem characterizes filters having the Baire property.

\begin{theorem}[Talagrand {[T1]}]\label{0.3}
For any filter ${\cal F}$ on $\omega$ the following conditions are equivalent:
\begin{enumerate}
\item ${\cal F}$ has the Baire property,
\item there exists a partition of
$\omega$ $\{I_{n} : n \in \omega \}$ into finite intervals such that
$\forall X \in {\cal F} \ \forall^{\infty}n \ X \cap I_{n} 
\neq \emptyset$.~$\QED$ \end{enumerate}
\end{theorem}   

Our first goal is to describe certain family
of $\mu_{\hat{p}}$-null sets which will be used
to cover $\mu_{\hat{p}}$-measurable filters.

{\sc Definition}
Set $H \subseteq 2^{\omega}$ is called {\em small}
with respect to the measure $\mu_{\hat{p}}$ if there exists a
partition of $\omega$ into pairwise
disjoint intervals $\{ I_{n} :n \in \omega\}$ and a sequence
$\{ J_{n} :n \in \omega\}$
such that
\begin{enumerate}
\item $J_{n} \subseteq 2^{I_{n}}$ for $n \in \omega$ ,
\item $H \subseteq \{x \in 2^{\omega}: \exists^{\infty}n \
x \rest I_{n} \in J_{n} \}$,
\item $\mu_{\hat{p}}(\{x \in 2^{\omega}: \exists^{\infty}n \ x \rest
I_{n}
\in J_{n}
\}) = 0$.
\end{enumerate}

Denote the set occuring in
$2)$ and $3)$ by $(I_{n},J_{n})_{n=1}^{\infty}$.
Notice that by Borel-Cantelli lemma we can replace condition (3) by the
equivalent :
$$3'.\hspace{0.1in} \sum_{n=1}^{\infty} \mu_{\hat{p}}(\{ x \in 2^{\omega} :
x \rest I_{n} \in J_{n} \}) < \infty .$$
The following generalizes the theorem from [Ba].
\begin{theorem}\label{0.4}
\begin{enumerate}
\item Every $\mu_{\hat{p}}$-null set is a union of two $\mu_{\hat{p}}$-small
sets,
\item There exists $\mu_{\hat{p}}$-null set which is not small.
\end{enumerate}
\end{theorem}
{\sc Proof}:
For completeness we sketch the proof of the first part.
Fix $\mu_{\hat{p}}$ for some
$\hat{p}=\{p_{n} : n \in \omega\}$ and let
$H \subseteq 2^{\omega}$ be a $\mu_{\hat{p}}$-null set.
The following claim is implicit in [O].
\begin{claim}
$\mu_{\hat{p}}(H)=0$  iff there exists a sequence
$\{F_{n} \subseteq 2^{n} : n \in \omega \}$ such that
\begin{enumerate}
\item $H \subseteq \{x \in 2^{\omega}: \exists^{\infty}n \ x \rest n
\in F_{n} \}  ,$
\item $\sum_{n=1}^{\infty} \mu_{\hat{p}}(\{x \in 2^{\omega}:
x \rest n \in  F_{n}\}) < \infty$.
\end{enumerate}
\end{claim}
{\sc Proof}:
The only difference between this and the definition of a small set is that
``domains" of $F_{n}$'s are not disjoint.

$\leftarrow$ This implication is an
immediate consequence of Borel-Cantelli lemma.

$\rightarrow$ Since $\mu_{\hat{p}}(H)=0$ there are
open sets $\{G_{n} :n \in \omega\}$ covering $H$ such
that $\mu_{\hat{p}}(G_{n}) < \frac{1}{2^{n}}$
for $n \in \omega$.
Write each $G_{n}$ as a union of disjoint basic sets i.e.
$$G_{n} = \bigcup_{m \in \omega}[s^{n}_{m}] \hbox{ for }  n \in \omega .$$
Let
$F_{n} = \{s \in 2^{n} : s=s_{k}^{l}$ for some $ k,l \in \omega \}$   for
$n \in \omega$.
Verification of $1)$ and $2)$ is straightforward.~$\QED$

Using the above claim and the assumption that
$\mu_{\hat{p}}(H)=0$ we can find a sequence
$\{F_{n} :n \in \omega \}$ satisfying the conditions
$1)$ and $2)$ of the claim. Fix a sequence of
positive reals $\{\varepsilon_{n} : n \in \omega \}$ such
that $\sum_{n=1}^{\infty} \varepsilon_{n} < \infty$  and let
$q(n) =
p_{1}^{-1}\cdot p_{2}^{-1}\cdot \ldots\cdot p_{n}^{-1}$ for
$n \in \omega$.

Define two sequences $\{n_{k},m_{k}:k \in \omega\}$ as follows:
$n_{0}=0$,
$$m_{k+1} = \min\{u > n_k:
q(n_{k})\cdot\sum_{j=u}^{\infty} \mu_{\hat{p}}(\{x \in 2^{\omega} :
x \rest j \in  F_{j}\}) < \varepsilon_{k}\}  ,$$
$$n_{k+1} = \min\{u >m_{k+1} : q(m_{k+1})\cdot\sum_{j=u}^{\infty}
\mu_{\hat{p}}(\{x \in 2^{\omega} : x \rest j \in F_{j}\}) <
\varepsilon_{k}\}  ,$$
and let
$$I_{k}=[n_{k},n_{k+1})  \hbox{ and }  I'_{k}=[m_{k},m_{k+1})
\hbox{ for }
k \in \omega .$$

Let
$$s \in  J_{k} \hbox{ iff }     s \in 2^{I_{k}} \hbox{ and }
\exists i \in [m_{k+1},n_{k+1})\  \exists t \in F_{i} \
s \rest \dom(t) \cap \dom(s)=$$
       $$ t \rest \dom(t) \cap \dom(s)  .$$
Similarly
$$s \in  J'_{k} \hbox{ iff } s \in 2^{I'_{k}} \hbox{ and } \exists i \in
[n_{k},m_{k+1})\
\exists t \in F_{i} \  s \rest \dom(t) \cap \dom(s)=$$
       $$ t \rest \dom(t) \cap \dom(s) . $$
It remains to show that $(I_{k},J_{k})_{k=1}^{\infty}$ and
$(I'_{k},J'_{k})_{k=1}^{\infty}$ are small with respect to
$\mu_{\hat{p}}$ and that their union  covers $H$.
Consider set $(I_{k},J_{k})_{k=1}^{\infty}$. Notice that
$$\mu_{\hat{p}}(\{x \in 2^{\omega} : x \rest I_{k} \in  J_{k}\}) \leq$$
$$\leq \mu_{\hat{p}}(\{x \in 2^{\omega} : \exists i \in [m_{k+1},n_{k+1}) \
\exists t \in F_{i} \  x \rest [m_{k+1},i)=
t \rest [m_{k+1},i)\}) \leq$$
$$q(n_{k})\cdot
\sum^{n_{k+1}-1}_{i=m_{k+1}} \mu_{\hat{p}}(\{x \in 2^{\omega} :
\exists t \in F_{i} \  x \rest i=t\}) \leq  \varepsilon_{k}.$$
Since $\sum_{n=1}^{\infty} \varepsilon_{n} < \infty$ this shows
that the set $(I_{n},J_{n})_{n=1}^{\infty}$ is $\mu_{\hat{p}}$-null.

Analogous argument works for the other set.
Finally we have that
$$H \subseteq (I_{n},J_{n})_{n=1}^{\infty} \cup (I'_{n},J'_{n})_{n=1}^{\infty}
.$$
To see this suppose that $x \in H$. Then
the set $X = \{n \in \omega: x \rest n \in F_{n}\}$ is infinite.
Thus either
$$X \cap \bigcup_{k=1}^{\infty} [m_{k+1},n_{k+1})\hbox{ is  infinite   or } $$
$$X \cap \bigcup_{k=1}^{\infty} [n_{k},m_{k+1})  \hbox{ is  infinite }.$$
Without loss of generality we can assume that it is the first case. But it
means that $x \in  (I_{n},J_{n})_{n=1}^{\infty}$ because if $x \rest n
\in  F_{n}$ and $n \in [m_{k+1},n_{k+1})$ then by the
definition of $J_{k}$ there is $t \in J_{k}$ such
that $x \rest [n_{k},n_{k+1})=t$. We are done since it
happens infinitely many times.~$\QED$

\section{Measurable filters}
In this section we characterize $\mu_{\hat{p}}$-measurable filters on $\omega$.

\begin{theorem}\label{1.1}
Let ${\cal F}$ be a filter on $\omega$. Then
${\cal F}$ is $\mu_{\hat{p}}$-measurable  iff ${\cal F}$ is
$\mu_{\hat{p}}$-small.
\end{theorem}
{\sc Proof}:
$\leftarrow$  This implication is obvious.

$\rightarrow$   Let ${\cal F}$ be a  $\mu_{\hat{p}}$-measurable filter.
Fix a sequence $\{\varepsilon_{n} : n \in \omega \}$
of positive reals such that $\sum_{n=1}^{\infty} 2^{n} \cdot \varepsilon_{n} <
\infty$.

For $J \subseteq 2^{<\omega}$ let $V(J) =
\{x \in  2^{\omega} : \exists s \in J \ s \subset x\}$.

By \ref{0.1} we know that ${\cal F}$ can be
covered by some $\mu_{\hat{p}}$-null set $H \subseteq 2^{\omega}$.
Applying \ref{0.4} we can find two $\mu_{\hat{p}}$-small sets
covering $H$. In fact as in the proof
of \ref{0.4} we can find sequences
$\{n_{k},m_{k} : k \in \omega\}$ and families
$\{J_{k}, J'_{k} : k \in \omega\}$ such that
\begin{enumerate}
\item $n_{k} < m_{k+1} < n_{k+1} < m_{k+2}$  for $k \in \omega$,
\item $J_{k} \subseteq 2^{[n_{k},n_{k+1})}  ,\
J'_{k} \subseteq 2^{[m_{k},m_{k+1})}$ for $k \in \omega$ ,
\item $\mu_{\hat{p}}(V(J_{k})) < \varepsilon_{k}$ ,
$\mu_{\hat{p}}(V(J'_{k})) < \varepsilon_{k}$  for $k \in \omega$ ,
\item $H \subseteq ([n_{k},n_{k+1}),J_{k})_{k=1}^{\infty}
\cup ([m_{k},m_{k+1}),J'_{k})_{k=1}^{\infty}$.
\end{enumerate}
Now for $k \in \omega$ define
$$S_{k} = \left\{t \in 2^{[n_{k},m_{k+1})} :
\mu_{\hat{p}}(V(\{s \in 2^{[m_{k+1},n_{k+1})} : t^{\frown}s \in  J_{k}\})) >
\frac{1}{2^{k}}\right\} .$$
Notice that for all $k \in \omega$
$$\frac{1}{2^{k}} \mu_{\hat{p}}(V(S_{k})) \leq  \mu_{\hat{p}}(J_{k}) <
\varepsilon_{k} \hbox{    hence}$$
$$\mu_{\hat{p}}(V(S_{k})) \leq 2^{k}\cdot\varepsilon_{k} .$$
Similarly if for $k>1$ we define
$$S'_{k} = \left\{t \in 2^{[n_{k},m_{k+1})} :
\mu_{\hat{p}}(V(\{s \in 2^{[m_{k},n_{k})} :
s^{\frown}t \in  J'_{k}\})) > \frac{1}{2^{k}} \right\}$$ then
$\mu_{\hat{p}}(V(S'_{k})) \leq  2^{k}\cdot\varepsilon_{k}$.

Thus we have three $\mu_{\hat{p}}$-small sets
$$H_{1} = ([n_{k},n_{k+1}),J_{k})_{k=1}^{\infty}  , $$
$$H_{2} = ([m_{k},m_{k+1}),J'_{k})_{k=1}^{\infty}\hbox{  and }$$
$$H_{3} = ([n_{k},m_{k}),S_{k} \cup S'_{k})_{k=1}^{\infty} .$$
If ${\cal F} \subseteq H_{2} \cup H_{3}$ we are done since it is
easy to see that it is a $\mu_{\hat{p}}$-small set. So
assume that there exists $X \in {\cal F}$ such
that $X \in  H_{1}$ but $X \not \in  H_{2} \cup H_{3}$.
Since $X \in  H_{1}$ we
have an infinite sequence $\{k_{u}:u \in \omega\}$ such that
$$\forall u \in \omega\   X \rest [n_{k_{u}},n_{k_{u}+1})
\in  J_{k_{u}} .$$
Define for $u \in \omega$
$$U_{u} = [m_{k_{u}+1},n_{k_{u}+1}) \hbox{  and}$$
$$T_{u} = \{s \in 2^{U_{u}} : X \rest [n_{k_{u}},m_{k_{u}+1})
^{\frown}s \in
J_{k_{u}} \  or\  s^{\frown}
X \rest [n_{k_{u}+1},m_{k{u}+2}) \in  J'_{k_{u}+1}\} .$$
We have to check that the set $(U_{u},T_{u})_{u=1}^{\infty}$
is $\mu_{\hat{p}}$-small. Consider sufficiently
large $u \in \omega$. Since
$X \rest [n_{k_{u}},n_{k_{u}+1}) \in  J_{k_{u}}$ and
$X \rest [n_{k_{u}},m_{k_{u}+1}) \not \in  S_{k_{u}} \cup
S'_{k_{u}}$
($u$ is large) we have $\mu_{\hat{p}}(V(T_{u})) < 2^{-u}$.
\begin{claim}
${\cal F}  \subseteq  (U_{u},T_{u})_{u=1}^{\infty}$.
\end{claim}
{\sc Proof}:
Suppose that ${\cal F}$ is not contained in this set and let
$Y \in {\cal F} -(U_{u},T_{u})_{u=1}^{\infty}$.
Define $Z \in  2^{\omega}$ as follows
$$Z(n) = \left\{ \begin{array}{ll}
X(n) & \hbox{if } n \in \bigcup_{u \in \omega} U_{u}\\
Y(n) & \hbox{ otherwise }
\end{array} \right. \hbox{ for } n \in \omega .$$
Notice that $Z \in {\cal F}$ since $X \cap Y \subseteq Z$.
We will show that $Z \not \in  H_{1} \cup H_{2}$ which gives
a contradiction. Consider an
interval $I_{m} = [n_{m},n_{m+1})$. If $m \neq k_{u}$ for
every $u \in \omega$ then $I_{m} \cap \bigcup_{u \in \omega}
U_{u} = \emptyset$ and $Z \rest I_{m}
\not \in  J_{m}$ since $Z \rest I_{m} =
X \rest I_{m}$ for such m's.
On the other hand if $m=k_{u}$ for some
$u \in \omega$ then
$X \rest I_{m} \in J_{m}$ but by the choice of $X \
Z \rest [n_{k_{u}},m_{k_{u}+1}) =
X \rest [n_{k_{u}},m_{k_{u}+1})$ has only
few extensions inside $J_{k_{u}}$ (since $X \not \in H_{3}$).
In fact if $Z \rest I_{m} \in J_{m}$
then $Z \rest [m_{k_{u}+1},n_{k_{u}+1})$
has to be an element of $T_{u}$. But this is
impossible since $Z \rest [m_{k_{u}+1},n_{k_{u}+1})
=Y \rest [m_{k_{u}+1},n_{k_{u}+1})
\not \in  T_{u}$ for sufficiently large
$u \in \omega$. Hence for all except finitely
many $m \in \omega$ $Z \rest I_{m} \not \in J_{m}$
which means that $Z \not \in H_{1}$.
Similarly, using the second clause in the definition of
$H_{3}$ we prove that $Z \not \in H_{2}$.
That finishes the proof since
$(U_{u},T_{u})_{u=1}^{\infty}$ is a $\mu_{\hat{p}}$-small set.~$\QED$

As a corollary we get:

\begin{theorem}\label{1.2}
For any filter ${\cal F}$ the following conditions are equivalent:
\begin{enumerate}
\item ${\cal F}$ is $\mu_{\hat{p}}$-measurable,
\item there exists a family $\{{\cal A}_{n} : n \in \omega \}$ such that
\begin{enumerate}
\item ${\cal A}_{n}$ consists of finitely many
finite subsets of $\omega$ for all $n \in \omega$,
\item $\bigcup{\cal A}_{n} \cap  \bigcup {\cal A}_{m} = \emptyset$
whenever $n \neq m$,
\item  $\sum_{n=1}^{\infty} \mu_{\hat{p}}(\{X \subseteq \omega :
\exists a \in {\cal A}_{n} \  a \subset X \}) < \infty$ ,
\item $\forall X \in {\cal F}\ \exists^{\infty}n\
\exists a \in {\cal A}_{n} \  a \subset X $.
\end{enumerate}
\end{enumerate}
\end{theorem}
{\sc Proof}:
$2) \rightarrow 1)$   This implication is obvious.

$1) \rightarrow 2)$. Assume that ${\cal F}$ is a measurable filter.
Then by the previous theorem ${\cal F} \subset (I_{n},J_{n})
_{n=1}^{\infty}$ for some $\mu_{\hat{p}}$-small set
$(I_{n},J_{n})_{n=1}^{\infty}$.
Define for $n \in \omega$
$$J'_{n} = \{s \in J_{n} : \forall u \in 2^{I_{n}} \
(s^{-1}(1) \subseteq u^{-1}(1) \rightarrow  u \in J_{n}\} .$$
\begin{claim}
${\cal F} \subseteq (I_{n},J'_{n})_{n=1}^{\infty}$.
\end{claim}
{\sc Proof}:
Suppose not. Let $X \in {\cal F} - (I_{n},J'_{n})_{n=1}^{\infty}$. It
is not very hard to see that there
exists a set $X' \supseteq X$ which
does not belong to $(I_{n},J_{n})_{n=1}^{\infty}$.
Contradiction.~$\QED$

Identify elements of $J'_{n}$ with subsets of $I_{n}$ and let
$${\cal A}_{n} = \{a \subseteq I_{n} : a \hbox{ is }
\subseteq-\hbox{minimal  element  of } J'_{n}\}\hbox{ for } n \in \omega.$$
Obviously ${\cal F} \subseteq  \{X \subseteq \omega : \exists^{\infty}n \
\exists a \in {\cal A}_{n} \  a \subset X\}$
and the family $\{{\cal A}_{n} : n \in \omega \}$ has properties $a)
-d)$.~$\QED$

If a family $\{{\cal A}_{n} : n \in \omega \}$ has   the properties
$a)-c)$ denote the set  $\{X \subseteq \omega : \exists^{\infty}n \
\exists a \in {\cal A}_{n} \
a \subset X\}$ by $({\cal A}_{n})_{n=1}^{\infty}$.
Characterization proved above can be interpreted as follows:
Filter  ${\cal F}$ is
$\mu_{\hat{p}}$-null iff there exists a
sequence of independent ``tests" $\{{\cal A}_{n} : n \in \omega\}$ such that
every element of ${\cal F}$ passes infinitely many of them.
Condition $(c)$ is a necessary
requirement for such a set to have measure zero.
Using \ref{1.2} one can prove that
\begin{theorem}\label{1.3}
Let ${\cal F}$ be a filter on $\omega$.
\begin{enumerate}
\item Let $\mu_{\hat{p}}$ and $\mu_{\hat{q}}$ be two
measures such that $p_{n} \leq q_{n}$ for all except finitely many
$n$.
Then $\mu_{\hat{p}}({\cal F})=0$ 
whenever $\mu_{\hat{q}}({\cal F})=0$.
\item Let $\mu_{\hat{p}}$ and $\mu_{\hat{q}}$ be two measures such
that ${\cal F}$ is nonmeasurable with respect to
both of them. Define $r_{n} = \min\{p_{n},q_{n}\}$ for
$n \in \omega$. Then ${\cal F}$ is $\mu_{\hat{r}}$-nonmeasurable.
\end{enumerate}
\end{theorem}
{\sc Proof}:
1) This can be showed be direct computation. A more sophisticated
but shorter is the following:

Suppose that $\hat{p}$ and $\hat{q}$ are two sequences such that
$p_{n} \leq q_{n}$ for $n \in \omega$. It is enough to show that
for any set $A = ({\cal A}_{n})_{n=1}^{\infty}$ , $\mu_{\hat{p}}(A) = 0$
whenever $\mu_{\hat{q}}(A)=0$.

Let $\varphi : P(\omega) \times P(\omega) \longrightarrow P(\omega)$
be the mapping defined as $\varphi(X,Y) = X \cup Y$ for $X,Y \in P(\omega)$.
Let
$$r_{n} = 1 - \frac{1-q_{n}}{1-p_{n}} \ for \ n \in \omega.$$

As in \ref{2.5} we show that
$\varphi^{-1}(\mu_{\hat{q}}) = \mu_{\hat{p}} \times \mu_{\hat{r}}$.

Since $\mu_{\hat{q}}(A)=0$ we have that $\mu_{\hat{p}} \times \mu_{\hat{r}}
(\varphi^{-1}(A))=0$. Therefore by Fubini's theorem there is
set $X \subseteq \omega$ such that
$$\mu_{\hat{p}}(\{Y \subseteq \omega: X \cup Y \in
({\cal A}_{n})_{n=1}^{\infty}
\})=0 .$$
But $\{Y \subseteq \omega: X \cup Y \in
({\cal A}_{n})_{n=1}^{\infty}
\} \supseteq ({\cal A}_{n})_{n=1}^{\infty}$.

2) Suppose that $X \subseteq \omega$ is an infinite set. Call $X$,
${\cal F}$-positive if the family ${\cal F} \cup \{X\}$ generates a
proper filter.

If $X$ is ${\cal F}$-positive let
$${\cal F}_{X} = \{ X \cap Y : Y \in {\cal F}\}$$
be a trace of ${\cal F}$ on $X$.
We will use the following fact:
\begin{claim}
For every filter ${\cal F}$ and measure $\mu_{\hat{p}}$
the following conditions are equivalent.
\begin{enumerate}
\item ${\cal F}$ is $\mu_{\hat{p}}$-nonmeasurable,
\item ${\cal F}_{X}$ is $\mu_{\hat{p} \rest X}$-nonmeasurable
for every ${\cal F}$-positive set $X \subseteq \omega$,
\item there exists ${\cal F}$-positive set $X \subseteq \omega$ such that
${\cal F}_{X}$ is $\mu_{\hat{p} \rest X}$-nonmeasurable.
\end{enumerate}
\end{claim}
{\sc Proof}:
1) $\rightarrow 2)$ Suppose that $\mu_{\hat{p} \rest X}
({\cal F}_{X})=0$ for some ${\cal F}$-positive set $X$.
Then by \ref{1.2} there exists $\mu_{\hat{p} \rest X}$-small
set $A = ({\cal A}_{n})_{n=1}^{\infty} \subset 2^{X}$ which
covers ${\cal F}_{X}$. Since $A$ is upwards-closed it is easy to see that
this set covers ${\cal F}$ as well and is $\mu_{\hat{p}}$-small.

$2) \rightarrow 3)$ Obvious.

$3) \rightarrow 1)$ Suppose that ${\cal F}_{X}$ is $\mu_{\hat{p} \rest
X}$-nonmeasurable for some $X \subseteq \omega$.
Notice that
$${\cal F}_{X} \times 2^{\omega - X} \subseteq {\cal F}.$$
That finishes the proof since $\mu_{\hat{p}} = \mu_{\hat{p} \rest X}
\times \mu_{\hat{p} \rest \omega - X}$.~$\QED$

Now we can finish the proof of \ref{1.3}. Suppose that ${\cal F}$ is
nonmeasurable with respect to measures $\mu_{\hat{p}}$ and $\mu_{\hat{q}}$.
Let $r_{n}= \min \{p_{n}, q_{n} \}$ for $n \in \omega$. We show that
${\cal F}$ is $\mu_{\hat{r}}$-nonmeasurable.
Let $X = \{ n \in \omega : p_{n} = r_{n}\}$. Clearly either $X$ or
$\omega - X$ is ${\cal F}$-positive. Without loss of generality we can assume
that we are in the first case. Applying the above claim and using the fact
that $\mu_{\hat{p} \rest X}$ and $\mu_{\hat{r} \rest X}$
are the same measures on $2^{X}$ we get the desired conclusion.~$\QED$

{\sc Definition}
Filter ${\cal F}$ is called {\em rapid} if for every increasing function
$f \in \omega^{\omega}$ there exists $X \in {\cal F}$ such that
$|X \cap f(n)| \leq n$ for $n \in \omega$.

As another application we get a simple proof of the following result
of Mokobodzki.
\begin{theorem}[Mokobodzki]\label{rapid}
Every rapid filter is Lebesque nonmeasurable.
\end{theorem}
{\sc Proof}:
Let ${\cal F}$ be a rapid filter. Suppose that ${\cal F}$
is covered by a set of form
$\{X \subset \omega: \exists^{\infty}n \ \exists a \in {\cal A}_{n}
\ a \subset X \}$ where $\{{\cal A}_{n}: n \in \omega\}$ is a family as
in \ref{1.2}.
Without loss of generality we can assume that
for all $n \in \omega,$
$$\mu(\{X \subseteq \omega: \exists a \in {\cal A}_{n} \ a \subset X\}) <
\frac{1}{2^{n+1}} $$
and that
$$\max\{\max(a) : a \in {\cal A}_{n}\} \geq
\min\{\min(a) : a \in {\cal A}_{m}\} \hbox{ for } n \geq m .$$
In particular it means that no set in ${\cal A}_{n}$ has less than
$n+1$ elements.
Define $f(n) =
\max\{\max(a) : a \in {\cal A}_{n}\} $ for $n \in \omega$ and let
$Z \in {\cal F}$ be such that
$|Z \cap f(n)|\leq n$ for all $n \in \omega$.
We immediately get that
$$ Z \not \in
\{X \subset \omega: \exists^{\infty}n \ \exists a \in {\cal A}_{n}
\ a \subset X \} .$$
Contradiction.~$\QED$

Before we go any further let us study the possible strenthening of \ref{1.2}.
For
simplicity we work with standard measure on $2^{\omega}$.
Suppose that $({\cal A}_{n})_{n=1}^{\infty}$  is a small set.
Notice that for given $n \in \omega$
$$\mu(\{X \subseteq \omega : \exists a \in {\cal A}_{n}\  a \subset X\}) \leq
\sum_{a \in {\cal A}_{n}} \frac{1}{2^{|a|}} - \sum_{a,b \in {\cal A}_{n}}
\frac{1}{2^{|a \cup b|}}+ \ldots .$$
Therefore it is natural to ask whether condition
$(c)$ in \ref{1.2} can be replaced by
the condition
$$\sum_{n=1}^{\infty} \sum_{a \in {\cal A}_{n}} \frac{1}{2^{|a|}} < \infty$$
or in general by
$$\sum_{n=1}^{\infty} \sum_{a \in {\cal A}_{n}} \prod_{i \in a} p_{i} <
\infty .$$
Surprisingly the answer turns out to be negative -- the following example was
found by M. Talagrand ([T3]).
\begin{theorem}[Talagrand]\label{1.4}
Assume CH. There exists a measurable filter
${\cal F}$ such that for every sequence $\{J_{n} :
 n \in \omega \}$ of finite subsets of $\omega$ satisfying
$$\sum_{n=1}^{\infty} \frac{1}{2^{|J_{n}|}} < \infty$$
there exists $X \in {\cal F}$ which contains no set
$J_{n}$ for $n \in \omega$.
\end{theorem}
{\sc Proof}:
Let us start with the following observation:
\begin{lemma}\label{1.5}
Let $I$ be a finite set of
size $2n$ for some $n \in \omega$ and let $\lambda_{I}$ be the counting
measure on the set $Z(I) = \{I' \subset I : |I'| = n \}$.
Suppose that $C$ is a subset of $I$. Then
$$\lambda_{I}(\{I'\in Z(I): C \subset I'\}) \leq \frac{1}{2^{|C|}} .$$
\end{lemma}
{\sc Proof}:
Suppose that $|C| = m$. Then the left-hand side is equal to
$$H(m)=\frac{\left( \begin{array}{c}
2n-m\\n-m \end{array} \right )}{\left( \begin{array}{c}
2n\\n \end{array} \right ) } \hbox{ so } H(m) = \frac{n-m}{2n-m}\cdot H(m-1)
\leq \frac{H(m-1)}{2} . \ \QED$$

\begin{lemma}\label{1.6}
Let $\{I_{n}:n \in \omega \}$ be a sequence of
pairwise disjoint subsets of $\omega$ each of them
having even number of elements.
Suppose that $\{J_{n} : n \in \omega \}$ is a
sequence of finite subsets of $\omega$ satisfying
$$\sum_{n=1}^{\infty} \frac{1}{2^{|J_{n}|}} < \infty.$$
Then for $n \in \omega$ there are
sets $I'_{n} \subset I_{n}$ of size $|I_{n}|/2$ such that
$\bigcup_{n \in \omega} I'_{n}$  contains no set $J_{n}$ for $n \in \omega$.
\end{lemma}
{\sc Proof}:
Provide $Z(I_{n})$ with the counting measure
$\lambda_{I_{n}}$ and ${\bold P} = \prod_{n \in \omega} Z(I_{n})$ with
the product 
measure $\lambda = \prod_{n \in \omega} \lambda_{I_{n}}$.
Using lemma \ref{1.5} we get that for
every $k \in \omega$
$$\lambda(\{\{I'_{n}: n \in \omega \} \in {\bold P} : J_{k}  \subseteq
\bigcup_{n \in \omega} I'_{n}\}) \leq
\prod_{n \in \omega} \lambda_{I_{n}}(\{I'_{n} \in Z(I_{n}):
J_{k} \cap I_{n} \subseteq I'_{n}\}) \leq$$
$$\leq \prod_{n \in \omega}
\frac{1}{2^{|I_{n} \cap J_{k}|}} = \frac{1}{2^{|J_{k}|}} .$$
Therefore
$$\lambda(\{\{I'_{n}:n \in \omega\} \in {\bold P}:
\forall k \in \omega\  J_{k} \not \subset \bigcup_{n \in \omega} I'_{n}\})
\geq \prod_{k \in \omega} (1- \frac{1}{2^{|J_{k}|}}) > 0 .$$
In particular the set of sequences we are looking
for has positive $\lambda$-measure.~$\QED$

{\sc Construction of the filter}

Let $\{I_{k,l} : k,l \in \omega\}$ be a family of pairwise disjoint sets
such that
$|I_{k,l}| = 2^{k}$ for $k,l \in \omega$.
Let $\{J_{n}^{\xi} : n \in \omega, \xi < \omega_{1}\}$ be an enumeration of
all sequences such that $$\sum_{n=1}^{\infty} \frac{1}{2^{|J_{n}^{\xi}|}} <
1 .$$ 
Construct by induction a sequence
$\{X_{\xi} : \xi<\omega_{1}\}$ of subsets of $\omega$ such that:
\begin{enumerate}
\item $J_{n}^{\xi} \not \subset X_{\xi}$  for $n \in \omega, \xi <
\omega_{1}$ , 
\item family $\{X_{\eta} : \eta<\xi\}$ has finite
intersection property for $\xi < \omega_{1}$ ,
\item for every $\xi<\omega_{1}$ and $\eta_{1}, \ldots, \eta_{n} < \xi$ there
exists a sequence of natural numbers
$\{a_{k} : k \in \omega\}$ such that
$\lim_{k \rightarrow \infty} a_{k} = \infty$ and $|X_{\eta_{1}}
\cap \cdots \cap X_{\eta_{n}} \cap I_{k,l}| \geq a_{k}$ for $l \in \omega$.
\end{enumerate}

Notice that it is enough to finish the proof:
let ${\cal F}$ be the filter generated by
the family $\{X_{\xi} : \xi < \omega_{1}\}$. Clearly
${\cal F}$ avoids every small set $({\cal A}_{n})_{n=1}^{\infty}$ such that
$\sum_{n=1}^{\infty} \sum_{a \in {\cal A}_{n}} 2^{-|a|} < \infty$.
Moreover
${\cal F}$ is null since ${\cal F}$ is contained in the set
$$\{X \subseteq \omega : \exists k \in \omega \ \forall l \in \omega
\ X \cap I_{k,l} \neq \emptyset\} $$
which is null.

Therefore assume that $\{X_{\beta} : \beta < \alpha <\omega_{1}\}$ are
already constructed. Order those sets in
order type $\omega$ say $\{Y'_{n}: n \in \omega\}$ and define
$$Y_{n} = Y'_{1} \cap \cdots \cap Y'_{n} \hbox{  for } n \in \omega .$$
By the induction hypothesis there are sequences
$\{a_{k}^{n}: k,n \in \omega \}$ such that
$\lim_{k \rightarrow \infty}
a_{k}^{n} = \infty$ for $n \in \omega$ and
$|Y_{n} \cap I_{k,l}| \geq a_{k}^{n}$ for $k,l,n \in \omega$.
Find a sequence $\{k_{n}: n \in \omega \}$ such that
$\lim_{n \rightarrow \infty} a_{k_{n}}^{n} = \infty$.
Let
$$X_{n} = Y_{n} \cap  \bigcup_{l \in \omega}
\bigcup_{j=k_{n}}^{k_{n+1}} I_{j,l} \hbox{ for } n \in \omega .$$
Let $X'_{\alpha} = \bigcup_{n \in \omega} X_{n}$. Now apply \ref{1.6} to
the sequence $\{J_{n}^{\alpha} : n \in \omega\}$ and partition
$\{X'_{\alpha} \cap I_{k,l} : k,l \in \omega\}$ to get a sequence
$\{I'_{k,l} : k,l \in \omega\}$.
Let
$$X_{\alpha} = \bigcup_{k,l \in \omega} I'_{k,l} .$$
Verification that $X_{\alpha}$ is the
element we are looking for is straightforward:
clearly $X_{\alpha}$ intersects every set $Y_{n}$
and avoids the sequence $\{J_{n}^{\alpha} : n \in \omega \}$.~$\QED$

\begin{theorem}
Every $\mu_{\hat{p}}$-measurable filter extends to a
$\mu_{\hat{p}}$-measurable filter 
which does not have the Baire property. 
\end{theorem}
{\sc Proof}: Let ${\cal F}$ be a measurable filter. By \ref{1.2}    we can find
a family $\{{\cal A}_{n} : n \in \omega\}$ such that
${\cal F} \subset ({\cal A}_{n})^{\infty}_{n=1}$.
For $X \subset \omega$ let $A_{X} = \{n \in \omega : \ \exists a \in {\cal
 A}_{n} \ a \subset X\}$.
It is easy to see that the family $\{A_{X} : X \in {\cal F}\}$ has
finite intersection property. Let ${\cal G}$ be any ultrafilter (or
filter which does not have the Baire property)
containing this family.
Define
$${\cal H} = \{X \subseteq \omega : A_{X} \in {\cal G}\} .$$
It is not very hard to see that the filter ${\cal H}$ has required
properties.~$\QED$ 

\section{Intersections of filters} This section is
devoted to the problem of measurability of the intersection
of family of filters. Let us start with countable case.
\begin{theorem}[Talagrand {[T1]}]\label{2.1}
Intersection of countably
many nonmeasurable filters is a nonmeasurable filter.~$\QED$
\end{theorem}
\begin{theorem}[Talagrand {[T1]}]\label{2.2}
Intersection of countable family of filters without the Baire property
is a filter 
without the Baire property.
Martin's Axiom implies that intersection of less
than $2^{\aleph_{0}}$ filters without the Baire
property does not have the Baire property.~$\QED$ 
\end{theorem}
Surprisingly the second part of the above theorem does not generalize when
category is replaced by measure. In fact we have the following:
\begin{theorem}[Fremlin {[F]}]\label{2.3}
Assume Martin's Axiom. Then there exists a family of Lebesgue nonmeasurable
filters of cardinality $2^{\aleph_{0}}$ such
that every uncountable subfamily has measurable
intersection.~$\QED$ 
\end{theorem}
The next theorem shows that the above pathology cannot happen if we assume
stronger measurability properties.
\begin{theorem}\label{2.4}
Assume Martin's Axiom. Let $\mu_{\hat{p}}$ be a
measure such that $\lim_{n \rightarrow \infty} p_{n}~=~0$ and let
$\{{\cal F}_{\xi} :
\xi < \lambda < 2^{\omega}\}$ be a family of $\mu_{\hat{p}}$-nonmeasurable
filters.
Then
$$\bigcap_{\xi < \lambda} {\cal F}_{\xi}
\hbox{ is a Lebesgue nonmeasurable filter.}$$
\end{theorem}
{\sc Proof}:
In fact we show that $\bigcap_{\xi<\lambda} {\cal F}_{\xi}$ is
$\mu_{\hat{q}}$-nonmeasurable
for any sequence $\hat{q}$ such that
$\lim_{n \rightarrow \infty} \frac{q_{n}}{p_{n}} = \infty$.
Let $\hat{q}=\{q_{n}:n \in \omega\}$ be a sequence satisfying the
above condition and let $({\cal A}_{n})_{n=1}^{\infty}$ be
any $\mu_{\hat{q}}$-small set. For given
$X \subseteq \omega$ let
$({\cal A}_{n}-X)_{n=1}^{\infty} = (\{a-X:a \in {\cal
A}_{n}\})_{n=1}^{\infty}$. 
Notice that if $X \in ({\cal A}_{n})_{n=1}^{\infty}$ then
$({\cal A}_{n} - X)_{n=1}^{\infty} = 2^{\omega}$.

Define by
induction sequences $\{X_{\xi}:\xi \leq \lambda\}$ and
$\{\hat{p}^{\xi}:\xi \leq \lambda\} \subset \Re^{\omega}$ such that
\begin{enumerate}
\item $X_{\xi} \in  {\cal F}_{\eta}$  for $\eta < \xi < \lambda$,
\item $X_{\xi} - X_{\eta}$ is finite for
$\xi < \eta \leq \lambda$ ,
\item $p_{n} = p_{n}^{\lambda} <p^{\xi}_{n} \leq
\frac{1}{2}p^{\eta}_{n} \leq q_{n}$ for
$\eta < \xi$
and all but finitely many $n \in \omega$,
\item $\lim_{n \rightarrow \infty} \frac{p^{\xi}_{n}}{p_{n}} = \infty$
for $\xi < \lambda$,
\item $\mu_{\hat{p}^{\xi}}(({\cal A}_{n}-X_{\xi})^{\infty}_{n=1}) = 0$ for
$\xi \leq \lambda$.
\end{enumerate}
It is easy to see that it is enough to finish
the proof: by $1)$ and $2)      \  X_{\lambda} \in
\bigcap_{\xi<\lambda} {\cal F}_{\xi}$
and $X_{\lambda} \not \in  ({\cal A}_{n})_{n=1}^{\infty}$ by $5)$ and the
remark above.

Suppose that $\{X_{\xi}: \xi < \alpha \}$ and $\{\hat{p}^{\xi}: \xi <
\alpha \}$ 
are already constructed and satisfy
conditions $1)-5)$.

{\sc Case 1}
$\alpha = \beta + 1$

Let $\varphi :2^{\omega} \times 2^{\omega} \longrightarrow  2^{\omega}$
be defined
as $\varphi(X,Y)(n)=\max\{X(n),Y(n)\}$. Notice that $\varphi(X,Y)$
is essentially the same as $X \cup Y$.

Define $p^{\alpha}_{n} = 1-\sqrt{1-p^{\beta}_{n}}$
for $n \in \omega$ and let $\nu$ be a measure on
$2^{\omega} \times 2^{\omega}$ defined as
$\mu_{\hat{p}^{\alpha}} \times \mu_{\hat{p}^{\alpha}}$.
As in \ref{2.5} we show that
$$\nu  = \varphi^{-1}(\mu_{\hat{p}^{\beta}}) .$$
Since by the induction hypothesis
$\mu_{\hat{p}^{\beta}}(({\cal A}_{n}-X_{\beta})_{n=1}^{\infty})=0$
we have that $\nu(\varphi^{-1}(({\cal A}_{n}-X_{\beta})
^{\infty}_{n=1}))=0$. We also have that ${\cal F}_{\beta}$ is
$\mu_{\hat{p}^{\alpha}}$-nonmeasurable because for almost every
$n \in \omega \   p_{n} \leq \frac{1}{2}p^{\beta}_{n}
\leq  1- \sqrt{1-p^{\beta}_{n}}$.
Therefore by Fubini theorem there exists
$X \in {\cal F}_{\beta}$ such that
$$\mu_{\hat{p}^{\alpha}}(\{Y \subseteq \omega : \varphi(X,Y)
\in ({\cal A}_{n}-X_{\beta})_{n=1}^{\infty}\})=0.$$
Let $X_{\alpha}=X_{\beta} \cup X$.

 By the above remarks
we have $\mu_{\hat{p}^{\alpha}}(({\cal A}_{n}-X_{\alpha})_{n=1}^{\infty}))=0$.
It is easy to check that other conditions are satisfied as well.

{\sc Case 2}
$\alpha$ is a limit ordinal.

For this case we use Martin's Axiom: first we construct a
sequence $\hat{p}^{\alpha}$
satisfying $3),4)$ and $5)$ and then $X_{\alpha}$ satisfying
$1)$ and $2)$.

Let ${\cal P}$ be the following notion of forcing:
$${\cal P} = \{\langle s,k,H \rangle :s \in {\bold Q}^{<\omega} , k
\in \omega , 
H \in [\alpha]^{<\omega}\ \& $$
$$\forall m>k \ \frac{\min\{p^{\xi}_{m} : \xi \in H\}}{p_{m}} > lh(s) \} .$$
({\bf Q} is the set of rationals).
For any $\langle s,k,H\rangle ,\langle s',k',H'\rangle  \in {\cal P}$ define
$$\langle s,k,H\rangle  \leq \langle s',k',H'\rangle \hbox{ iff }  s'
\subseteq s 
\ \& \  k' \leq k \ \& \ H' \subseteq H\  \&\
\hbox{ for  all } n \geq k'$$
$$lh(s')\cdot p_{n} \leq s(n) \leq \min\{p^{\xi}_{n} :
\xi \in H'\} .$$
One easily checks that ${\cal P}$ is ccc.
Define for $\xi < \alpha \  D_{\xi} = \{\langle s,k,H\rangle : \xi \in H\}$.
It is easy to see that these sets are dense
in ${\cal P}$. If ${\bold G}$ is a filter which intersects all of them define
$$\hat{p}' = \bigcup \{s:\langle s,k,H\rangle  \in {\bold G}\} \hbox{  and }
p^{\alpha}_{n} = \frac{p'_{n}}{2}\hbox{  for } n \in \omega .$$
It is not very hard to check that this is a sequence we were looking for.
Now we construct $X_{\alpha}$.

By the induction hypothesis we know
that $\mu_{\hat{p}^{\xi}}(({\cal A}_{n}-X_{\xi})_{n=1}^{\infty}))=0$ for
$\xi < \alpha$. Therefore
by \ref{1.3} $\mu_{\hat{p}^{\alpha}}(({\cal
A}_{n}-X_{\alpha})_{n=1}^{\infty}))=0$
for $\xi < \alpha$ which is equivalent to
$\sum_{n=1}^{\infty} \mu_{\hat{p}^{\alpha}}(\{X \subset \omega :
\exists a \in {\cal A}_{n}-X_{\xi}\  a \subset X\}) < \infty$ for
$\xi < \alpha$.

Using Martin's Axiom find $\sum_{n=1}^{\infty} \varepsilon_{n} < \infty$ such
that
$$\forall \xi < \alpha \  \forall^{\infty}n\
\mu_{\hat{p}^{\alpha}}(\{X \subset \omega: \exists a \in {\cal A}_{n}
- X_{\xi}\ 
a \subset X\}) < \varepsilon_{n}.$$
Let ${\cal Q}$ be the following notion of forcing:
$${\cal Q} = \{\{\langle n_{\alpha_{1}},X_{\alpha_{1}}\rangle
,\langle n_{\alpha_{2}},X_{\alpha_{2}}\rangle,  \ldots,
\langle n_{\alpha_{k}},X_{\alpha_{k}}\rangle \}:
k \in \omega \ \&\  \alpha_{i} < \alpha\hbox{  for }    i<k \ \&$$
$$\hbox{ for    all } n \in \omega \
\mu_{\hat{p}^{\alpha}}(\{X \subset \omega :
\exists a \in {\cal A}_{n}-(\bigcup_{i<k}\langle
X_{\alpha_{i}}-n_{\alpha_{i}}\rangle
\ a \subset X\}) < \varepsilon_{n}\})\} .$$
For $p,q \in {\cal Q}$ define $p \leq q$ iff $p \supseteq q$.
\begin{claim}
${\cal Q}$ is ccc.
\end{claim}
{\sc Proof}:
Let $W \subseteq {\cal Q}$ be an uncountable family.
By ``thinning out" we can assume that there
are $k,n_{1},\ldots n_{k} \in \omega$ such that every
element of $W$ is of the form
$\{\langle n_{1},X_{\alpha_{1}}\rangle , \ldots, \langle
n_{k},X_{\alpha_{k}}\rangle \}$.
Observe that for every $X_{\alpha_{1}}, \ldots ,X_{\alpha_{j}}$
there is $n \in \omega$ such that
$\{\langle n,X_{\alpha_{1}}\rangle ,\ldots ,\langle n,X_{\alpha_{j}}\rangle \}
\in
{\cal Q}$.
This is because sets $\{X_{\beta}:  \beta < \alpha\}$ form an
increasing sequence.

Since $W$ is an uncountable antichain
we can find a number $n \in \omega$
and two conditions
$\{\langle n_{1},X_{\alpha_{1}}\rangle ,\ldots ,\langle
n_{k},X_{\alpha_{k}}\rangle \}\in W$ and
$\{\langle n_{1},X_{\beta_{1}}\rangle ,\ldots ,\langle n_{k},\
 X_{\beta_{k}}\rangle
\}
\in W$ such that
$\{\langle n,X_{\alpha_{1}}\rangle ,\ldots ,\langle n,X_{\alpha_{k}}\rangle
,\langle n,X_{\beta_{1}}\rangle ,\ldots, \langle n,X_{\beta_{k}}\rangle \} \in
{\cal
Q}$
and
$$n \cap \bigcup_{i<k} \langle X_{\alpha_{i}}-n_{\alpha_{i}}\rangle  =
n \cap \bigcup_{i<k} \langle X_{\beta_{i}}-n_{\beta_{i}}\rangle  .$$
Thus conditions
$\{\langle n_{1},X_{\alpha_{1}}\rangle ,\ldots, \langle
n_{k},X_{\alpha_{k}}\rangle \}$
and $\{\langle n_{1},X_{\beta_{1}}\rangle ,\ldots\-,\langle
n_{k},X_{\beta_{k}}\rangle \}$
are
compatible, which finishes the proof.~$\QED$

Let $D_{\xi}=\{p \in {\cal Q}: \exists n \in \omega\
\langle n,X_{\xi}\rangle \in p\}$ for $\xi < \alpha$.
It is easy to see that all sets  $D_{\xi}$ are
dense in ${\cal Q}$. Let ${\bold G}$ be a filter
intersecting all  $D_{\xi}$'s. Define
$$X_{\alpha} = \bigcup \{X_{\xi}-n_{\xi}: \exists p \in {\bold G} \
\langle n_{\xi},X_{\xi} \rangle \in p\}.$$
Verification that $X_{\alpha}$ satisfies conditions
$1)-5)$ is straightforward.~$\QED$

Notice that if the family of filters is countable we do not need Martin's
Axiom. Using the same method we can prove Talagrand's theorem from
the first section.
\begin{corollary}[Talagrand]
Let $\{{\cal F}_{n} :n \in \omega \}$ be a countable
family of $\mu_{\hat{p}}$-nonmeasurable filters.
Then
$$\bigcap_{n \in \omega} {\cal F}_{n} 
\hbox{ is  a }  \mu_{\hat{p}}\hbox{-nonmeasurable       filter.}$$
\end{corollary}
{\sc Proof}:
Let $({\cal A}_{n})_{n=1}^{\infty}$ be any $\mu_{\hat{p}}$-small set.
Construct a sequence $\{X_{n}:n \in \omega \}$ as in the proof
of 2.4 for measures
$\mu_{\hat{p}_{m}} \  m \in \omega$ where $p^{m}_{n} =
2^{-m}p_{n} $ for $n,m \in \omega$.
Set $X_{\omega}$ will witness
that $\bigcap_{n \in \omega} {\cal F}_{n}$ is not covered
by $({\cal A}_{n})_{n=1}^{\infty}$.
Use the fact that ${\cal F}_{n}$ is
$\mu_{\hat{p}^{m}}$-nonmeasurable for $n,m \in \omega$ [T1].~$\QED$

Therefore, if we have Martins's Axiom countable case generalizes
to uncountable
provided we have little bit stronger measurability hypothesis.

\section{Filters which are both null and meager}
This section is devoted to filters which are both null and meager.
Let ${\cal F}$ be a $\mu_{\hat{p}}$-measurable filter.
By \ref{1.2} ${\cal F}$ can be covered by some $\mu_{\hat{p}}$-small set
$({\cal A}_{n})_{n=1}^{\infty}$.
For $X \in \omega$ define $\supp(X)=\{n \in \omega
:\exists a \in {\cal A}_{n}\  a \subset X\}$ and let
$$ {\cal F}^{\star} = \{ \supp(X): X \in {\cal F}\}.$$
Notice that the definition of ${\cal F}^{\star}$ makes sense only in presence
of some covering $({\cal A}_{n})_{n=1}^{\infty}$ of ${\cal F}$.
It is easy to see that ${\cal F}^{\star}$ is a filter which
is a continuous image of  ${\cal F}$.

\begin{lemma}\label{3.1}
If ${\cal F}$ is $\mu_{\hat{p}}$-measurable and ${\cal F}^{\star}$ has
the Baire property then ${\cal F}$ can be covered by a
$\mu_{\hat{p}}$-null set of type $F_{\sigma}$.
\end{lemma}
{\sc Proof}:
Suppose that ${\cal F} \subseteq ({\cal A}_{n})_{n=1}^{\infty}$.
If ${\cal F}^{\star}$ has
the Baire property then using theorem \ref{0.3} we
can find a partition of $\omega \ \{I_{n}:n \in \omega \}$ such that
$$\forall X \in {\cal F}^{\star}\  \forall^{\infty}n\
X \cap I_{n} \neq  \emptyset .$$
As a consequence we get
$${\cal F} \subseteq  \bigcup_{n \in \omega}\
\bigcap_{m \geq n}\  \bigcup_{k \in I_{m}}\  \{X \subseteq \omega
:\exists a \in {\cal A}_{k} \  a \subset X\}.$$
The above set is a $\mu_{\hat{p}}$-null set of type $F_{\sigma}$.~$\QED$

\begin{corollary}
Every Borel (analytic) filter can be covered by a null set of
type $F_{\sigma}$.~$\QED$
\end{corollary}

Notice that if ${\cal F}$ can be covered
by a null set of type $F_{\sigma}$ then ${\cal F}$ is measurable
and has the
Baire property. The next theorem shows that the converse does not hold.
For simplicity let us work with Lebesgue measure.
\begin{theorem}
Assume that there exists a nonmeasurable filter having the Baire property.
Then there exists a filter which is both null and meager but cannot be covered
by a null set of type $F_{\sigma}$.
\end{theorem}
{\sc Proof}:
Let us first notice that the existence of nonmeasurable filter having the Baire
property follows from Martin's Axiom but it is not provable
in ZFC. (see [BGJS]).

Let ${\cal G}$ be a nonmeasurable filter
with the Baire property and ${\cal H}$ any filter without
the Baire property.

Let $\{I_{n}:n \in \omega \}$ be a partition witnessing
that ${\cal G}$ is meager. We can assume that
$|I_{n}|>n$ for $n \in \omega$. Define
$${\cal F} = \{X \in {\cal G} : \{n \in \omega :I_{n} \subset X\} \in
{\cal H}\} .$$
It is very easy to verify that ${\cal F}$ is a filter,
${\cal F}$ is meager since ${\cal F} \subseteq {\cal G}$ and
${\cal F}$ is
null since ${\cal F} \subseteq \{X \in \omega : \exists^{\infty}n\
I_{n} \subset X\}$ which is null.
We will show that ${\cal F}$ cannot be covered
by a null set of type $F_{\sigma}$. Let $K \subset 2^{\omega}$ be
such a set. First find an increasing sequence of closed sets
$\{C_{n} : n \in \omega \}$ such that
$K \subset \bigcup_{n \in \omega}C_{n}$ and $\mu(C_{n})=0$ for
$n \in \omega$.
Now for $n,m \in \omega$ define
$$C^{n}_{m} =\{s \in 2^{m} : [s] \cap C_{n} \neq \emptyset \} .$$
Let $\{k_{n}:n \in \omega \}$ be any sequence of natural numbers such that
$$\sum_{n=1}^{\infty} 2^{k_{n}}\cdot\mu(V(C^{n}_{k_{n+1}})) < \infty .$$
Define for $n \in \omega$
$$U_{n} = [k_{n},k_{n+1})\hbox{  and }$$
$$T_{n} = \{s \in 2^{U_{n}} : \exists t \in C^{n}_{k_{n+1}}\
s \rest U_{n} = t \rest U_{n}\}.$$
From the above definitions easily follows that
$$K \subseteq \{x \in 2^{\omega} : \forall^{\infty}n \
x \rest U_{n} \in T_{n}\} \subseteq (U_{n},T_{n})_{n=1}^{\infty}$$
and that the set $(U_{n},T_{n})_{n=1}^{\infty}$ is small.
Without loss of generality we can also assume that
$$\forall n \in \omega\ \exists m \in \omega \  I_{n} \subset U_{m} .$$
Since ${\cal G}$ is a nonmeasurable filter we
can find $X \in  {\cal G}-(U_{n},T_{n})_{n=1}^{\infty}$. Using theorem
\ref{0.3} and the fact that ${\cal H}$ does not have
the Baire property we can also find an
element $Y \in {\cal H}$ such that for some infinite set $S \subseteq \omega$
$$(\star) \hspace{0.1in} \bigcup_{n \in Y}
I_{n} \cap \bigcup_{n \in S} U_{n} = \emptyset.$$
Let
$$Z = X \cup  \bigcup_{n \in Y} I_{n}. $$
We will show that $Z \in  {\cal F} - K$ which finishes
the proof since $K$ is an arbitrary $F_{\sigma}$
set. Obviously $Z \in {\cal F}$ and $Z \not \in
\{x \in 2^{\omega} : \forall^{\infty}n\  x \rest U_{n} \in T_{n}\}$
because of $(\star)$ and the definition
of $X$.~$\QED$

\medskip

\begin{center}
{\it Received September 18, 1991}
\end{center}

\end{document}